\documentclass[draft,12pt]{elsarticle}

\makeatletter
\def\ps@pprintTitle{%
	\let\@oddhead\@empty
	\let\@evenhead\@empty
	\def\@oddfoot{}%
	\let\@evenfoot\@oddfoot}
\makeatother 

\usepackage{a4wide}

\usepackage{amssymb}
\usepackage{amsmath}
\usepackage{graphicx} 

\newtheorem{theorem}{Theorem}

\newtheorem{corollary}{Corollary}

\newtheorem{definition}{Definition}

\newtheorem{remark}{Remark}

\def\al{\alpha}

\def\({\left(}
\def\){\right)}
\def\[{\left[}
\def\]{\right]}

\newcommand\bin[2]{\left(\!\!\!\begin{array}{c}#1\\#2\end{array}\!\!\!\right)}

\def\N{\mathbb{N}}

\def\C{\mathbb{C}}

\date{}

\newcommand\beq{  \begin{equation}}
\newcommand\eeq{\end{equation}  }

\numberwithin{equation}{section}

\begin{document}

\begin{frontmatter}

\title{Recurrence relations for   Apostol-Bernoulli , -Euler and -Genocchi polynomials of higher order }
\author{Marc Pr\'evost \textsuperscript {1,2,3}}

\address{
\begin{trivlist}
\item\textsuperscript {1} Univ Lille Nord de France, F-59000 Lille,
France

\item\textsuperscript {2} Univ. Littoral C\^ote d’Opale, EA 2597 - LMPA - Laboratoire de
Math\'ematiques Pures et
Appliqu\'ees Joseph Liouville, F-62228 Calais, France

\item\textsuperscript {3} CNRS,FR 2956, France
FR 2956
\end{trivlist}}

\ead{marc.prevost@univ-littoral.fr}

\begin{abstract}
In \cite{luo2006,luosri2005}, Luo and Srivastava introduced some generalizations of the Apostol -Bernoulli polynomials and   the Apostol-Euler polynomials. The main object of this paper is to extend the result of \cite{prevost2010} to these generalized polynomials. More precisely, using the Pad\'{e} approximation of the exponential
function, we obtain recurrence relations for    Apostol-Bernoulli, Euler   and also Genocchi polynomials of higher order. As an application we prove lacunary relation for some particular cases.
 
\end{abstract}


\end{frontmatter}

\section{Introduction, definition and notations}

 The generalized Bernoulli polynomials $B_n^{(\alpha)}(x)$ of order $\al \in \C$, the generalized Euler polynomials 
$E_n^{(\alpha)}(x)$ of order $\al \in \C$ and  the generalized Genocchi  polynomials 
$G_n^{(\alpha)}(x)$ of order $\al \in \C$, each of degree $n$, are respectively defined by the folloing generating functions (see  \cite{erdelyi1953}, vol. III, p.253 and seq., \cite{luke1969}, Section 2.8 and \cite{luosri2011}:

\begin{eqnarray}
\left(\frac{t}{ \,e^{t}-1}\right)^\alpha e^{x\,t} &=&\sum_{k=0}^{\infty }\,{B}%
_{k}^{(\alpha)}(x)\frac{t^{k}}{k!}, \;\;\;\;	(\left\vert t
\right\vert <2\pi ; 1^\al:=1, \al \in \C),
\end{eqnarray}
\begin{eqnarray}
\left(\frac{2}{ \,e^{t}+1}\right)^\alpha e^{x\,t} &=&\sum_{k=0}^{\infty }\,{E}%
_{k}^{(\alpha)}(x)\frac{t^{k}}{k!}, \;\;\;\;	(\left\vert t
\right\vert <\pi ; 1^\al:=1, \al \in \C),
\end{eqnarray}

and 
\begin{eqnarray}
\left(\frac{2t}{ \,e^{t}+1}\right)^\alpha e^{x\,t} &=&\sum_{k=0}^{\infty }\,{G}%
_{k}^{(\alpha)}(x)\frac{t^{k}}{k!}, \;\;\;\;	(\left\vert t
\right\vert <\pi ; 1^\al:=1,\al \in \C),
\end{eqnarray}

The classical Bernoulli polynomials $B_n(x)$, the classical Euler polynomials $E_n(x)$  and   the classical Genocchi polynomials are given by

$B_n(x):=B_n^{(1)}(x), E_n(x):=E_n^{(1)}(x),G_n(x):=G_n^{(1)}(x),(n\in \N)$ respectively.

Moreover, the Bernoulli $B_n$'s numbers, Euler $E_n$'s numbers and Genocchi $G_n$'s numbers   are given by:
$B_n:=B_n(0), E_n:=2^n E_n(\frac{1}{2}), G_n:=G_n(0).$

These polynomials and numbers play a fundamental role in various branches of mathematics including combinatorics, number theory and special functions.

Q.M.  Luo and Srivastava  introduced the Apostol-Bernoulli polynomials of higher order (also called   generalized Apostol-Bernoulli polynomials):

\begin{definition}(Luo and Srivastava, \cite{luosri2005})
The Apostol-Bernoulli polynomials $\mathcal{B}^{(\alpha)}%
_{k}(x;\lambda )$ of order $\alpha$ in the variable $x$ are defined by means of the generating function:
\begin{eqnarray}
\left(\frac{t}{\lambda \,e^{t}-1}\right)^\alpha e^{x\,t} &=&\sum_{k=0}^{\infty }\,\mathcal{B}%
_{k}^{(\alpha)}(x;\lambda )\frac{t^{k}}{k!}, \label{genapostbern}
\end{eqnarray}

$$(\left\vert t
\right\vert <2\pi { \it \;  when\; } \lambda=1; \left\vert t
\right\vert <\left|\log \lambda\right| , 1^\al :=1).$$
\end{definition}

According to the definition, by setting $\alpha=1$, we obtain the Apostol-Bernoulli polynomials ${\mathcal B}_k(x;\lambda)$. Moreover, we call ${\mathcal B}_k(\lambda):={\mathcal B}_k(0;\lambda)$ the Apostol-Bernoulli numbers.

Explicit representation of $\mathcal{B}_{k}^{(\alpha)}(x;\lambda )$ in terms of a generalization of the Hurwitz-Lerch zeta function  can be found in \cite{gargjainsri2006}.

In \cite{luo2006}, Luo introduced the Apostol-Euler polynomials of higher order $\alpha$.

\begin{definition}
 The Apostol-Euler   polynomials $\mathcal{E}^{(\alpha)}
_{k}(x;\lambda )$ of order real or complex $\alpha$ 
in the variable $x$ are
defined by means of the following generating function:%
\begin{eqnarray}
\left( \frac{2}{\lambda \,e^{t}+1}\right)^\alpha e^{x\,t} &=&\sum_{k=0}^{\infty }\mathcal{E}^{(\alpha)}
_{k}(x;\lambda )\,\frac{t^{k}}{k!}\;\;(  \left\vert t
\right\vert <\left|\log(-\lambda)\right|, 1^\al:=1).  \label{genapostoleuler}
\end{eqnarray}
\end{definition}

The Apostol-Euler polynomials ${\mathcal E}_k(x;\lambda)$ are given by ${\mathcal E}_k(x;\lambda):={\mathcal E}_k^{(1)}(x;\lambda)$. The Apostol-Euler numbers ${\mathcal E}_k( \lambda)$ are given by ${\mathcal E}_k( \lambda):=2^k{\mathcal E}_k( \frac{1}{2};\lambda)$.

 Some relations between  Apostol-Bernoulli and Apostol-Euler polynomials  of order $\alpha$ can be found in \cite{luosri2006}. For more results on these polynomials, the readers are referred to \cite{choiandsri2008,luo2009bis,luo2009}.

 In \cite{luo2009fourier} Luo  introduced and investigated the Apostol-Genocchi polynomials of order $\al$, which are defined as follows.
  
  \begin{definition}
  	The Apostol-Genocchi   polynomials $\mathcal{G}^{(\alpha)}
  	_{k}(x;\lambda )$ of order   $\alpha $
  	in the variable $x$ are
  	defined by means of the following generating function:%
  	\begin{eqnarray}
  	\left( \frac{2t}{\lambda \,e^{t}+1}\right)^\alpha e^{x\,t} &=&\sum_{k=0}^{\infty }\mathcal{G}^{(\alpha)}
  	_{k}(x;\lambda )\,\frac{t^{k}}{k!}\;\;\;(  \left\vert t
  	\right\vert <\left|\log(-\lambda)\right|, 1^\al:=1). \label{genapostolgenocchi}
  	\end{eqnarray}
  \end{definition}
  
  The Apostol-Genocchi polynomials ${\mathcal G}_k(x;\lambda)$ are given by ${\mathcal G}_k(x;\lambda):={\mathcal G}_k^{(1)}(x;\lambda)$. The Apostol-Genocchi numbers ${\mathcal G}_k( \lambda)$ are given by ${\mathcal G}_k( \lambda):={\mathcal G}_k( 0;\lambda)$.
  
  When $\lambda=1$ in (\ref{genapostbern}) and when $\lambda=-1$ in (\ref{genapostolgenocchi}), the order $\al$ of the generalized Apostol- Bernoulli polynomial 
  $\mathcal{B}_k^{(\al)}(x;\lambda)$ and the order $\al$ of the generalized Apostol- Genocchi polynomial 
  $\mathcal{G}_k^{(\al)}(x;\lambda)$ should tacitly be restricted to non negative values.

\bigskip

In this paper, we   consider Apostol-type polynomials of order $\al$, 
$\mathcal{ B}_k^{(\al)}(x;\lambda)$  
$\mathcal{ E}_k^{(\al)}(x;\lambda)$ and $\mathcal{ G}_k^{(\al)}(x;\lambda)$.

The aim of this paper is to apply  Pad\'{e} approximation to 
$e^{t}\ $in the generating functions (\ref{genapostbern}) ,  (\ref{genapostoleuler}) and (\ref{genapostolgenocchi}) to get  relations between Apostol-type polynomials of order
$\al$ depending on three parameters $n,\;m,\;p$ where $n$\ and $m$\ are respectively the degree of
the numerator and the degree of the denominator of the Pad\'{e} approximant
used to approximate the function $e^{t}$\ and $p$\ is some positive integer.

The paper is organized as follows. In the  next section, we recall the
definition of Pad\'{e} approximant to a general series and  its
expression for the case of the exponential function. In Section 3, we apply
Pad\'{e} approximation to prove the   recurrence relations (Theorem \ref{teo1} and \ref{teo2})
for the Apostol-Bernoulli, Apostol-Euler and Apostol-Genocchi polynomials of higher order.

\section{Pad\'{e} approximant}

In this section, we recall the definition of Pad\'{e} approximation to
general series and their expression in the case of the exponential function.
Given a function $f$\ with a Taylor expansion 
\begin{equation}
f(t)=\sum_{i=0}^{\infty }c_{i}t^{i}  \label{seriesf}
\end{equation}%
in a neighborhood of the origin, a Pad\'{e} approximant denoted $[n,m]_{f}$\
to $f$\ is a rational fraction of degree   $n$ (resp. $m$) for the  numerator (resp. the denominator):
\begin{equation*}
\lbrack n,m]_{f}(t)=\frac{\alpha _{0}+\alpha _{1}t+\cdots +\alpha _{n}t^{n}}{%
\beta _{0}+\beta _{1}t+\cdots +\beta _{m}t^{m}},
\end{equation*}%
whose Taylor expansion agrees with (\ref{seriesf}) as far as possible:

\begin{equation*}
\sum_{i=0}^{\infty }c_{i}t^{i}-\frac{\alpha _{0}+\alpha _{1}t+\cdots +\alpha
_{n}t^{n}}{\beta _{0}+\beta _{1}t+\cdots +\beta _{m}t^{m}}=O(t^{m+n+1}).
\end{equation*}%
In the general case, the resulting linear system has  unique solutions $%
\alpha _{i},\beta _{i}$ (see, e.g., \cite{bgm1996}).

  Pad\'e approximation is related with  convergence  acceleration \cite{brez2000,prevost1994,prevost1998,prevost2000},   continued fractions \cite{brez1976, jonesthron1980},   orthogonal polynomials,   quadrature formulas \cite{gautschi1985} and number theory.
 Moreover the denominators of Pad\'e approximants satisfy a three terms recurrence  \cite{brez1980, brez1990} and this property allows  finding  another proof of the irrationality of $\zeta(2)$ and $\zeta(3)$ \cite{prevost1996}.

\bigskip
\bigskip

\bigskip

  If $f(t)=e^{t}$\ then   
\begin{align*}
\lbrack n,m]_{f}(t):& =\frac{P^{(n,m)}(t)}{Q^{(n,m)}(t)}=\frac{%
_{1}F_{1}(-n,-m-n,t)}{_{1}F_{1}(-m,-m-n,-t)} \\
& =\frac{\sum_{j=0}^{n}\frac{(-n)_{j}}{(-m-n)_{j}}\frac{t^{j}}{j!}}{%
\sum_{j=0}^{m}\frac{(-m)_{j}}{(-m-n)_{j}}\frac{(-t)^{j}}{j!}}=\frac{%
\sum_{j=0}^{n}\left(\!\!\!
\begin{array}{c}
n \\ 
j%
\end{array}%
\!\!\!\right) (n+m-j)!t^{j}}{\sum_{j=0}^{m}\left(\!\!\!
\begin{array}{c}
m \\ 
j%
\end{array}%
\!\!\!\right) (n+m-j)!(-t)^{j}},
\end{align*}

\noindent where the Pochhammer symbol $(a)_{j}$\ is defined as 
\begin{align*}
(a)_{j}& =a(a+1)\cdots (a+j-1){{\mathrm{\;if\;}}}j\geq 1, \\
& =1{{\mathrm{\;if\;}}}j=0 \end{align*}
 and the hypergeometric series $_{1}F_{1}(a,b,z)$ is defined as ${_{1}F_{1}(a,b,z)}:=\sum_{k=0}^\infty \frac{(a)_k}{(b)_k}\frac{z^k}{k!}.$
In the sequel we write  $[n,m](t)$ for the Pad\'e approximant to $e^t$.
 The remainder term is defined by 
\begin{equation*}
R^{(n,m)}(t):=e^{t}-[n,m](t)=e^{t}-\frac{P^{(n,m)}(t)}{Q^{(n,m)}(t)}
\end{equation*}%
and satisfies

\begin{align*}
R^{(n,m)}(t)& =t^{m+n+1}\frac{e^{t}}{Q^{(n,m)}(t)}%
\int_{0}^{1}(x-1)^{m}x^{n}e^{-x\;t}dx \\
& =\frac{(-1)^{m}}{Q^{(n,m)}(t)}\sum_{j=0}^{\infty }\frac{t^{j+m+n+1}}{%
(j+m+n+1)\left(\!\!\! 
\begin{array}{c}
m+n+j \\ 
n%
\end{array}%
\!\!\!\right) j!}\;\;\; \\
& =\frac{(-1)^{m}\;n!}{Q^{(n,m)}(t)}\sum_{j=0}^{\infty }\frac{(m+j)!}{%
j!\,(m+n+j+1)!}t^{j+m+n+1}\; \\
& =O(t^{n+m+1}).
\end{align*}%
\noindent In the sequel, let

\begin{align*}
\alpha _{j}^{(n,m)}& :=\left(\!\!\! 
\begin{array}{c}
n \\ 
j%
\end{array}%
\!\!\!\right) {(n+m-j)!},\;\;\;0\leq j\leq n, \\
\beta _{j}^{(n,m)}& :=\left(\!\!\! 
\begin{array}{c}
m \\ 
j%
\end{array}%
\!\!\!\right) {(n+m-j)!}(-1)^{j},\;\;\;0\leq j\leq m, \\
\gamma _{j}^{(n,m)}& :={(-1)^{m}}\frac{n!\;(m+j)!}{(m+n+1+j)!\;j!},\;\;\;j\geq 0.
\end{align*}

\noindent  The polynomials $P^{(n,m)},Q^{(n,m)}$ and the product $R^{(n,m)}\;Q^{(n,m)}$ are then given by 
\begin{align*}
P^{(n,m)}(t)& =\sum_{j=0}^{n}\alpha _{j}^{(n,m)}t^{j}, \\
Q^{(n,m)}(t)& =\sum_{j=0}^{m}\beta _{j}^{(n,m)}t^{j}, \\
R^{(n,m)}(t)\;Q^{(n,m)}(t)& =\sum_{j=0}^{\infty }\gamma
_{j}^{(n,m)}t^{j+m+n+1}.
\end{align*}%

\section{\ Recurrence relations for Apostol-type  polynomials}

Let us recall a classical  method to derive a basic formula for Apostol-Bernoulli polynomials of order $\al$.

From the definition, 
\begin{eqnarray}
\left(\frac{t}{\lambda \,e^{t}-1}\right)^\alpha e^{(x+y)\,t} &=&\sum_{k=0}^{\infty }\,\mathcal{B}%
_{k}^{(\alpha)}(x+y;\lambda )\frac{t^{k}}{k!}, 
\end{eqnarray}
we write the left handside member of the previous equation as
\begin{eqnarray}
\left(\frac{t}{\lambda \,e^{t}-1}\right)^\alpha e^{(x+y)\,t} &=& 
\left(\frac{t}{\lambda \,e^{t}-1}\right)^\alpha e^{x\,t} e^{y \;t} 
\end{eqnarray}
and     substitute $e^{y \;t}$ by its Taylor expansion around $0$.
It arises

\beq
\left( \sum_{k=0}^{\infty }\,\mathcal{B}
_{k}^{(\alpha)}(x;\lambda )\frac{t^{k}}{k!}\right) \left( \sum_{k=0}^{\infty }\,   \frac{y^k\;t^{k}}{k!}\right)=\sum_{k=0}^{\infty }\,\mathcal{B}%
_{k}^{(\alpha)}(x+y;\lambda )\frac{t^{k}}{k!},
\eeq

 By identification, the following formula is proved:
 
 \beq
 \mathcal{B}
 _{n}^{(\alpha)}(x+y;\lambda )=\sum_{k=0}^{\infty }\,\bin{n}{k}\mathcal{B}%
 _{n-k}^{(\alpha)}(x;\lambda )y^k.\label{formule addition}
\eeq
In this section, the main idea is to replace the exponential function $e^{y\;t}$ in the generating functions by some Pad\'e approximant of Apostol-type polynomials (\ref{genapostbern}) ,
(\ref{genapostoleuler}) and (\ref{genapostolgenocchi}),  $e^{y\;t}$\ not by its Taylor
expansion around $0$\ but by its Pad\'{e} approximant $[n,m]$. This gives
the following main result.

\begin{theorem}
\label{teo1}

For all integers $m\geq 0,n\geq 0$, for arbitrary  real or complex parameter $\lambda$,
the Apostol-type polynomials ${\Lambda}_k^{(\al)}(x;\lambda) $of order $
\al$ satisfy

for $0\leq p\leq m+n,$
\beq
\sum_{k=0}^{n }\,
\dfrac{(n+m-k)!}{(p-k)!}\bin{n}{k}y^k{\Lambda}%
_{p-k}^{(\alpha)}(x;\lambda )
=\sum_{k=0}^{m  }\,  (-y)^k
\dfrac{(n+m-k)!}{(p-k)!}\bin{m}{k}{\Lambda}%
_{p-k}^{(\alpha)}(x+y;\lambda )\label{formule1teo1}
\eeq

for $ p\geq m+n+1,$
\begin{multline}
\sum_{k=0}^{n }\,
\dfrac{(n+m-k)!}{(p-k)!}\bin{n}{k}y^k{{\Lambda}}%
_{p-k}^{(\alpha)}(x;\lambda )
=\sum_{k=0}^{m }\,  (\!-\!y\!)^k
\dfrac{(n+m-k)!}{(p-k)!}\bin{m}{k}{\Lambda}
_{p-k}^{(\alpha)}(x+y;\lambda )-\\
\!\!\!\!\!\!\!\!\!\!\!\!\!\!\!\!\!\!\!(\!-\!1\!)^m \dfrac{n!m!}{p!}\sum_{k=0}^{\!p-\!m-\!n-\!1 }\,  \!y^{p-k}
\bin{p}{k}\bin{\!p\!-\!k\!-\!n\!-\!1}{m}{\Lambda}%
_{k}^{(\alpha)}(x;\lambda )\label{formule2teo1}
\end{multline}
where the ${\Lambda}%
_{k}^{(\alpha)}(x;\lambda )$'s are  the Apostol-Bernoulli, -Euler or -Genocchi polynomials of order $\al$.

\end{theorem}

We use the convention ${{\Lambda}}%
_{k}^{(\alpha)}(x;\lambda )=0$ for $k\leq -1$.

\bigskip

Particular cases:

1) If $p=m+n$, then 
\beq
\sum_{k=0}^{n }\, 
\bin{n}{k}y^k{\Lambda}%
_{p-k}^{(\alpha)}(x;\lambda )
=\sum_{k=0}^{m }\,  (-y)^k
\bin{m}{k}{\Lambda}%
_{p-k}^{(\alpha)}(x+y;\lambda ).
\eeq

2) If $m=0$, then 
\beq
{\Lambda}
_{p}^{(\alpha)}(x+y;\lambda )=\sum_{k=0}^{p}\bin{p}{k }y^k{\Lambda}
_{p-k}^{(\alpha)}(x;\lambda ).
\eeq

Formula (\ref{formule addition}) is then recovered.

3) If $n=0$, \beq
{\Lambda}
_{p}^{(\alpha)}(x;\lambda )=\sum_{k=0}^{p}\bin{p}{k }(-y)^k{\Lambda}
_{p-k}^{(\alpha)}(x+y;\lambda ).
\eeq
which is the dual  formula of the previous one.

\noindent \textbf{Proof of Theorem \ref{teo1}} 

We display the proof only for the case of Apostol-Bernoulli polynomials  of order $\al$.
Let us start from the generating functions (\ref{genapostbern}), (\ref{genapostoleuler}),(\ref{genapostolgenocchi})
$$\Phi(t,\lambda,\al)e^{x\;t}=\sum_{k=0}^\infty\Lambda_k^{(\al)}(x;\lambda)\dfrac{t^k}{k!}$$
where 
\begin{eqnarray*}
\Phi(t,\lambda,\al)&=&\left(\dfrac{t}{\lambda e^t-1}\right)^\al, \rm{  \;for\; Apostol-Bernoulli\; polynomials\; of \;order\; } \al,\\
\Phi(t,\lambda,\al)&=&\left(\dfrac{2}{\lambda e^t+1}\right)^\al, \rm{  \;for\; Apostol-Euler\; polynomials\; of \;order\; } \al,\\
\Phi(t,\lambda,\al)&=&\left(\dfrac{2t}{\lambda e^t+1}\right)^\al, \rm{  \;for\; Apostol-Genocchi\; polynomials\; of \;order\; } \al,
\end{eqnarray*}

and for some integers $n,m$, we replace $e^{y\;t}$ by its Pad\'e approximant previously defined:
\begin{eqnarray*}
e^{y\;t}&=&[n,m](y \;t)+R^{(n,m)}(y\;t)\\
&=&\dfrac{P^{(n,m)}(y\;t)}{Q^{(n,m)}(y\;t)}+R^{(n,m)}(y\;t).
\end{eqnarray*}

 We get $$\Phi(t,\lambda,\al)e^{x\;t}e^{y\;t}=\Phi(t,\lambda,\al)e^{x\;t}\left(\dfrac{P^{(n,m)}(y\;t)}{Q^{(n,m)}(y\;t)}+R^{(n,m)}(y\;t)\right)=\sum_{k=0}^\infty\Lambda_k^{(\al)}(x+y;\lambda)\dfrac{t^k}{k!}$$
 with
 
  \noindent $P^{(n,m)}(y\;t)=\displaystyle \sum_{k=0}^n \alpha _{k}^{(n,m)}y^k t^k=\sum_{k=0}^n\bin{n}{k}(m+n-k)!\dfrac{t^k}{k!}$,

  \noindent $Q^{(n,m)}(y\;t)=\displaystyle\sum_{k=0}^m \beta _{k}^{(n,m)}y^k t^k=\sum_{k=0}^m\bin{m}{k}(-1)^k (m+n-k)!\dfrac{t^k}{k!}$,

 \noindent and  
 
 \noindent $R^{(n,m)}(y\;t)\;Q^{(n,m)}(y\;t)=\displaystyle\sum_{k=0}^\infty \gamma_{k}^{(n,m)}(y\;t)^{m+n+1+k} =\sum_{k=0}^\infty (-1)^{m}\frac{n!\;(m+
 k)!}{(m+n+1+k)!\;k!}(y\;t)^{m+n+1+k}. $

This leads to 

 $$\sum_{k=0}^\infty\Lambda_k^{(\al)}(x;\lambda)\dfrac{t^k}{k!}\times
 \left(P^{(n,m)}(y\;t)+R^{(n,m)}(y\;t)\;Q^{(n,m)}(y\;t)\right) =
 Q^{(n,m)}(y\;t)\sum_{k=0}^\infty\Lambda_k^{(\al)}(x+y;\lambda)\dfrac{t^k}{k!}
 $$
  
 which gives, applying the Cauchy product for series, the following identity:
 
 $$\sum_{p=0}^\infty t^p \!\!\sum_{\substack{k+j=p \\ k \geq 0, j\geq 0} }\!\!
 \alpha _{j}^{(n,m)}\dfrac{\Lambda_k^{(\al)}(x;\lambda)}{k!}=
 \!\!\sum_{p=0}^\infty t^p  \!\!\sum_{\substack{k+j=p \\ k \geq 0, j\geq 0} }
 \!\!\beta _{j}^{(n,m)}\dfrac{\Lambda_k^{(\al)}(x+y;\lambda)}{k!}-
 \sum_{p=0}^\infty t^p  \!\!\sum_{\substack{k+j=p\\-m-n-1 \\ k \geq 0, j\geq 0} }\!\!
 \gamma _{j}^{(n,m)}\dfrac{\Lambda_k^{(\al)}(x+y;\lambda)}{k!}
 $$
 
 Comparing the coefficient of $t^p$ on both sides of the previous equation, we obtain the assertions  (\ref{formule1teo1}) and (\ref{formule2teo1})  of Theorem \ref{teo1}.

 \begin{flushright}
 	 $\blacksquare $

 \end{flushright}

 Theorem \ref{teo1} provides a recurrence 
 relation of length $\max (n,m)$\ from $\mathcal{B}^{(\al)}%
_{p-\max (n,m)}(x;\lambda )$ (resp. $\mathcal{E}^{(\al)}_{p-\max (n,m)}(x;\lambda
)$, $\mathcal{G}^{(\al)}_{p-\max (n,m)}(x;\lambda)$)  to $\mathcal{B}^{(\al)}_{p}(x;\lambda )$\ (resp. $\mathcal{E}^{(\al)}_{p}(x;\lambda )$,
$\mathcal{G}^{(\al)}_{p}(x;\lambda ))$
\ for $p$\ less than $m+n.$
 On the other hand, if $p$\ is greater than $m+n,$ we obtain also a
recurrence relation for the same Apostol-type polynomials of higher order, but with supplementary
first terms from $\mathcal{B}^{(\al)}_{0}(x;\lambda )$\ (resp. $\mathcal{E}^{(\al)}%
_{0}(x;\lambda )$, $\mathcal{G}^{(\al)}%
_{0}(x;\lambda )$) to $\mathcal{B}^{(\al)}_{p-m-n-1}(x;\lambda )$\ (resp. $\mathcal{E%
}^{(\al)}_{p-m-n-1}(x;\lambda )$, resp. $\mathcal{G%
}^{(\al)}_{p-m-n-1}(x;\lambda )$).

\bigskip Of course, Theorem \ref{teo1} is valid for classical Bernoulli, Euler and Genocchi  polynomials and also for Bernoulli, Euler and Genocchi numbers which are particular
Apostol-type polynomials. By this mean we recover  known results in the literature as shown in \cite{prevost2010}.

\bigskip

\bigskip
\bigskip

By using the same method  as in the proof of  Theorem \ref{teo1}, we establish now a recurrence formula for ${\Lambda}_{k}^{(\al+1)}(x;\lambda )$ and $\Lambda_{k}^{(\al)}(x;\lambda )$.
\begin{theorem}\label{teo2}
	
	For $n,m, p  \in \N$,
	
if $0\leq p\leq m+n,$ then the following relations for Apostol-type polynomials of consecutive order $\al$ and $\al+1$,
\begin{eqnarray}
	\sum_{j=0}^{\max(m,n)} (\lambda \al_j^{(n,m)}-\beta_j^{(n,m)})\dfrac{\mathcal{B}%
	_{p-j}^{(\alpha+1)}(x;\lambda )}{(p-j)!}
&=&\sum_{j=0}^{m} \beta_j^{(n,m)}\dfrac{\mathcal{B}%
	_{p-j-1}^{(\alpha)}(x;\lambda )}{(p-j-1)!},\label{formule1teo2}
\\
	\sum_{j=0}^{\max(m,n)} (\lambda \al_j^{(n,m)}+\beta_j^{(n,m)})\dfrac{\mathcal{E}%
	_{p-j}^{(\alpha+1)}(x;\lambda )}{(p-j)!}
&=&2\sum_{j=0}^{m} \beta_j^{(n,m)}\dfrac{\mathcal{E}%
	_{p-j}^{(\alpha)}(x;\lambda )}{(p-j)!},\label{formule2teo2}\\
	\sum_{j=0}^{\max(m,n)} (\lambda \al_j^{(n,m)}+\beta_j^{(n,m)})\dfrac{\mathcal{G}%
	_{p-j}^{(\alpha+1)}(x;\lambda )}{(p-j)!}
&=&2\sum_{j=0}^{m} \beta_j^{(n,m)}\dfrac{\mathcal{G}%
	_{p-j-1}^{(\alpha)}(x;\lambda )}{(p-j-1)!},\label{formule3teo2}
\end{eqnarray}
hold, 

where $\al_j^{(n,m)}=\bin{n}{j}(n+m-j)!$ and $\beta_j^{(n,m)}=(-1)^j \bin{m}{j}(n+m-j)!$.

For $  p\geq m+n+1,$

$$
\sum_{j=0}^{\max(m,n)} (\lambda \al_j^{(n,m)}-\beta_j^{(n,m)})\dfrac{\mathcal{B}%
	_{p-j}^{(\alpha+1)}(x;\lambda )}{(p-j)!}
=\sum_{j=0}^{m} \beta_j^{(n,m)}\dfrac{\mathcal{B}%
	_{p-j-1}^{(\alpha)}(x;\lambda )}{(p-j-1)!}-
\lambda \sum_{j=0}^{p-m-n-1} \gamma_j^{(n,m)}\dfrac{\mathcal{B}%
	_{p-m-n-j-1}^{(\alpha+1)}(x;\lambda )}{(p-m-n-j-1)!}
$$

$$
\sum_{j=0}^{\max(m,n)} (\lambda \al_j^{(n,m)}+\beta_j^{(n,m)})\dfrac{\mathcal{E}%
	_{p-j}^{(\alpha+1)}(x;\lambda )}{(p-j)!}
=2\sum_{j=0}^{m} \beta_j^{(n,m)}\dfrac{\mathcal{E}%
	_{p-j}^{(\alpha)}(x;\lambda )}{(p-j)!}-
\lambda \sum_{j=0}^{p-m-n-1} \gamma_j^{(n,m)}\dfrac{\mathcal{E}%
	_{p-m-n-j-1}^{(\alpha+1)}(x;\lambda )}{(p-m-n-j-1)!}
$$

$$
\sum_{j=0}^{\max(m,n)} (\lambda \al_j^{(n,m)}+\beta_j^{(n,m)})\dfrac{\mathcal{G}%
	_{p-j}^{(\alpha+1)}(x;\lambda )}{(p-j)!}
=2\sum_{j=0}^{m} \beta_j^{(n,m)}\dfrac{\mathcal{G}%
	_{p-j-1}^{(\alpha)}(x;\lambda )}{(p-j-1)!}-
\lambda \sum_{j=0}^{p-m-n-1} \gamma_j^{(n,m)}\dfrac{\mathcal{G}%
	_{p-m-n-j-1}^{(\alpha+1)}(x;\lambda )}{(p-m-n-j-1)!}
$$
where $\gamma_j^{(n,m)}=(-1)^m\dfrac{n!(m+j)!}{(m+n+1+j)!j!}$
\end{theorem}

\noindent \textbf{Proof}  

We only display the proof for the Apostol-Bernoulli polynomials of order $\al$.

We start from the generating function (\ref{genapostbern})
$$
\left(\frac{t}{\lambda \,e^{t}-1}\right)^\alpha e^{x\,t} =\sum_{k=0}^{\infty }\,\mathcal{B}%
_{k}^{(\alpha)}(x;\lambda )\frac{t^{k}}{k!}.$$
We multiply the two hand-side members of this equation by $\dfrac{t}{\lambda e^t-1}$ to get
$$\sum_{k=0}^{\infty }\,\mathcal{B}%
_{k}^{(\alpha+1)}(x;\lambda )\frac{t^{k}}{k!}=
\left(\frac{t}{\lambda \,e^{t}-1}\right) \sum_{k=0}^{\infty }\,\mathcal{B}%
_{k}^{(\alpha)}(x;\lambda )\frac{t^{k}}{k!}$$
and 
$$(\lambda \,e^{t}-1)\sum_{k=0}^{\infty }\,\mathcal{B}%
_{k}^{(\alpha+1)}(x;\lambda )\frac{t^{k}}{k!}=
  \sum_{k=0}^{\infty }\,\mathcal{B}%
_{k}^{(\alpha)}(x;\lambda )\frac{t^{k+1}}{k!},$$
in which we replace $e^t$ by its Pad\'e approximant.
This leads to $$
(\lambda P^{(n,m)}(t)  -Q^{(n,m)}(t)+\lambda R^{(n,m)}(t) Q^{(n,m)}(t))
\sum_{k=0}^{\infty }\,\mathcal{B}%
_{k}^{(\alpha+1)}(x;\lambda )\frac{t^{k}}{k!}=Q^{(n,m)}(t)\sum_{k=0}^{\infty }\,\mathcal{B}%
_{k}^{(\alpha)}(x;\lambda )\frac{t^{k+1}}{k!},$$

where 

 $P^{(n,m)}(t)=\displaystyle\sum_{k=0}^n \alpha _{k}^{(n,m)}t^k$,
 $Q^{(n,m)}(t)=\displaystyle\sum_{k=0}^m \beta _{k}^{(n,m)}t^k$ and  
  $R^{(n,m)}(t)\;Q^{(n,m)}(y\;t)=\displaystyle\sum_{k=0}^\infty \gamma_{k}^{(n,m)}(t)^{m+n+1+k} .$
  
  Theorem \ref{teo2} concerned with Apostol-Bernoulli polynomials od order $\al$ is obtained by equating the coefficient of $t^p$..

 \begin{flushright}
 $\blacksquare $

 \end{flushright}

\section{Applications}

\bigskip\ In this section, we will consider the value of the parameter $p\,\ 
$with respect to $n,m.$

Let us first consider the particular values $p=m+n$ and $\al=0$. Then we can prove the following corollary.

\begin{corollary}
	\label{coro0}
For $m\geq 0,n\geq 0,$%
\begin{eqnarray}
	\lambda \sum_{k=0}^{n}\bin{n}{k}{B}_{m+k}(x;\lambda
	)\!-\!(\!-\!1)^{m}\sum_{k=0}^{m}(-1)^{k}\bin{m}{k}{B}_{n+k}(x;\lambda )
\!\!\!\!	&=&\!\!\!\!((m+n)x\!-\!n)x^{n-1}(x\!-\!1)^{m-1},\label{formule1 cor1} \\
	 \lambda \sum_{k=0}^{n}\bin{n}{k}{E}_{m+k}(x;\lambda
	)+(-1)^{m}\sum_{k=0}^{m}(-1)^{k}\bin{m}{k}{E}_{n+k}(x;\lambda )
	\!\!\!\!&=&\!\!\!\!2x^{n}(x-1)^{m},\label{formule2 cor1}\\
	\lambda \sum_{k=0}^{n}\bin{n}{k}{G}_{m+k}(x;\lambda
	)\!+\!(\!-\!1)^{m}\sum_{k=0}^{m}(\!-\!1)^{k}\!
	\bin{m}{k}{G}_{n+k}(x;\lambda )
	\!\!\!\!&=&\!\!\!\!2((m\!+\!n)x\!-\!n)x^{n\!-\!1}(x\!-\!1)^{m\!-\!1}.\label{formule3 cor1}
\end{eqnarray}

\end{corollary}

{\bf  Proof}

The two  formulas (\ref{formule1 cor1}), (\ref{formule2 cor1}) have been proved in \cite{prevost2010}.
We prove the formula  (\ref{formule3 cor1}).

If $p=m+n$ and $\al=0$ in (\ref{formule3teo2}), then 

$$
 \sum_{j=0}^{\max(m,n)} (\lambda \al_j^{(n,m)}+\beta_j^{(n,m)})\dfrac{{G}%
	_{m+n-j} (x;\lambda )}{(m+n-j)!}
=2\sum_{j=0}^{m} \beta_j^{(n,m)}\dfrac{x^
	{m+n-j-1} }{(m+n-j-1)!}
$$
since $\mathcal{G}%
_{k}^{(1)}(x;\lambda )={G}%
_{k}(x;\lambda )$ and $\mathcal{G}%
_{k}^{(0)}(x;\lambda )=x^k.$

Thus replacing $\al_j^{(n,m)}$ , $\beta_j^{(n,m)}$ by their expression, we obtain the following relation:

\begin{eqnarray*}
\lambda \sum_{j=0}^{n}  \bin{n}{j}G_{m+n-j}(x;\lambda)
+\sum_{j=0}^{m}  (-1)^j \bin{m}{j}G_{m+n-j}(x;\lambda)
\!\!\!\!&=&\!\!\!\!2\sum_{j=0}^{m}(-1)^j \bin{m}{j}(m+n-j) x^{m+n-j-1} \\
&=&2((m+n)x-n)x^{n-1}(x-1)^{m-1}
\end{eqnarray*}

 \bigskip 
 \begin{flushright}
 	$\blacksquare $
 \end{flushright}

  Let $s\geq 1 $. Setting   $p=m+n-s$   in 
      Theorem \ref{teo2},   we will immediately get   Corollary \ref{coro2}.

\begin{corollary}\label{coro2} 

For $n\in \mathbb{N}$, $m\in \mathbb{N}$, $s\in \mathbb{N}$ such that $\;1\leq s\leq n+m,$

\noindent $\lambda \displaystyle\sum_{k=0}^{n}
\bin{n}{k}
 (m\!-\!s\!+\!k\!+\!1)_s\mathcal{B}^{(\al+1)}_{m-\!s+\!k}(x;\lambda )\!-\!(-1)^{m}\!\displaystyle\sum_{k=0}^{m}\!(-1)^{k}\bin{m}{k} \!(n\!-\!s\!+\!k\!+\!1)_{s}\mathcal{B}^{(\al+1)}_{n-s+k}(x;\lambda )=$
 \begin{flushright}
 $\!(-1)^{m}\!\displaystyle\sum_{k=0}^{m}\!(-1)^{k}\bin{m}{k} \!(n\!-\!s\!+\!k\!)_{s+1}\mathcal{B}^{(\al)}_{n-s+k-1}(x;\lambda ),$
 \end{flushright}

\noindent $\lambda \displaystyle\sum_{k=0}^{n}\bin{n}{k} (m\!-\!s\!+\!k\!+\!1)_{s}\mathcal{E}^{(\al+1)}_{m-s+k}(x;\lambda )\!+\!(-1)^{m}\!\sum_{k=0}^{m}(-1)^{k}\bin{m}{k} (n-s+k+1)_{s}\mathcal{E}^{(\al+1)}_{n-s+k}(x;\lambda )=$
\begin{flushright}
	$\!2(-1)^{m}\!\displaystyle\sum_{k=0}^{m}\!(-1)^{k}\bin{m}{k} \!(n\!-\!s\!+\!k\!+1)_{s}\mathcal{E}^{(\al)}_{n-s+k}(x;\lambda ),$
\end{flushright}
 
 \noindent $\lambda \displaystyle\sum_{k=0}^{n}
 \bin{n}{k}
 (m\!-\!s\!+\!k\!+\!1)_s\mathcal{G}^{(\al+1)}_{m-\!s+\!k}(x;\lambda )\!+\!(-1)^{m}\!\sum_{k=0}^{m}\!(-1)^{k}\bin{m}{k} \!(n\!-\!s\!+\!k\!+\!1)_{s}\mathcal{G}^{(\al+1)}_{n-s+k}(x;\lambda )=$
 \begin{flushright}
 	$2(-1)^{m}\!\displaystyle\sum_{k=0}^{m}\!(-1)^{k}\bin{m}{k} \!(n\!-\!s\!+\!k\!)_{s+1}\mathcal{G}^{(\al)}_{n-s+k-1}(x;\lambda ),$
 \end{flushright}
 where the Pochhammer symbol $(a)_{j}$\ is defined as 
 \begin{align*}
 (a)_{j}& =a(a+1)\cdots (a+j-1){{\mathrm{\;if\;}}}j\geq 1, \\
 & =1{{\mathrm{\;if\;}}}j=0 \end{align*}
\end{corollary}

\begin{remark}
	
For $\lambda =1$ and $s=1$, Corollary  \ref{coro2}  reduces to the following formulas,

 $\displaystyle\sum_{k=0}^{n} \bin{n}{k} (m+k){B}^{(\al+1)}_{m+k-1}(x )\!-\!(-1)^{m}\!\sum_{k=0}^{m}\!(-1)^{k} \bin{m}{k} \!(n+k){B}^{(\al+1)}_{n+k-1}(x)=$\\
 \begin{flushright}
 $\!(-1)^{m}\!\displaystyle\sum_{k=0}^{m}\!(-1)^{k}\bin{m}
 {k} \!(n+k-1)(n+k){B}^{(\al)}_{n+k-2}(x),$
 \end{flushright}

$ \displaystyle\sum_{k=0}^{n}\bin{n}{k} (m+k){E}^{(\al+1)}_{m+k-1}(x )\!+\!(-1)^{m}\!\displaystyle\sum_{k=0}^{m}(-1)^{k}\bin{m}
{k} (n+k){E}^{(\al+1)}_{n+k-1}(x)=$
 \begin{flushright}
 $\displaystyle 2(-1)^{m}\sum_{k=0}^{m}(-1)^k\bin{m}{k}(n+k)
 {E}^{(\al)}_{n+k-1}(x)$,
\end{flushright}

$\displaystyle\sum_{k=0}^{n} \bin{n}{k} (m+k){G}^{(\al+1)}_{m+k-1}(x )\!+\!(-1)^{m}\!\sum_{k=0}^{m}\!(-1)^{k} \bin{m}{k} \!(n+k){G}^{(\al+1)}_{n+k-1}(x)=$
\begin{flushright}
	$\displaystyle2(-1)^{m}\!\sum_{k=0}^{m}\!(-1)^{k}\bin{m}
	{k} \!(n+k-1)(n+k){G}^{(\al)}_{n+k-2}(x ),$
\end{flushright}

 which extend   Kaneko's formula \cite{kaneko1995}  to Bernoulli, Euler and Genocchi polynomials of higher order $\al$.

\end{remark}

\begin{remark}
	
	If $p=m+n+r$, with $r\geq 1$, similar relations exist. We only display the case $r=1$:
	
	$$
\lambda	\sum_{j=0}^{n} \bin{n}{j}\dfrac{\mathcal{B}%
		_{m+j+1}^{(\alpha+1)}(x;\lambda )}{m+j+1}+(-1)^{m+1} 
	\sum_{j=0}^{m}(-1)^j \bin{m}{j}\dfrac{\mathcal{B}%
		_{n+j+1}^{(\alpha+1)}(x;\lambda )}{n+j+1}$$$$
	=(-1)^m \sum_{j=0}^{m} (-1)^j\bin{m}{j}\mathcal{B}%
		_{n+j}^{(\alpha)}(x;\lambda )-
	\lambda (-1)^m \dfrac{n!m!}{(n+m+1)!}\mathcal{B}%
	_{0}^{(\alpha+1)}(x;\lambda ) 
	$$

	$$
	\lambda	\sum_{j=0}^{n} \bin{n}{j}\dfrac{\mathcal{E}%
		_{m+j+1}^{(\alpha+1)}(x;\lambda )}{m+j+1}+(-1)^m
	\sum_{j=0}^{m} (-1)^j \bin{m}{j}\dfrac{\mathcal{E}%
		_{n+j+1}^{(\alpha+1)}(x;\lambda )}{n+j+1}$$$$
	=2(-1)^m \sum_{j=0}^{m} (-1)^j \bin{m}{j}\dfrac{\mathcal{E}%
		_{n+j+1}^{(\alpha)}(x;\lambda )}{n+j+1}-
	\lambda (-1)^m \dfrac{n!m!}{(n+m+1)!}\mathcal{E}%
	_{0}^{(\alpha+1)}(x;\lambda ) 
	$$

		$$
		\lambda	\sum_{j=0}^{n} \bin{n}{j}\dfrac{\mathcal{G}%
			_{m+j+1}^{(\alpha+1)}(x;\lambda )}{m+j+1}+(-1)^{m} 
		\sum_{j=0}^{m}(-1)^j \bin{m}{j}\dfrac{\mathcal{G}%
			_{n+j+1}^{(\alpha+1)}(x;\lambda )}{n+j+1}$$$$
		=2(-1)^m \sum_{j=0}^{m} (-1)^j\bin{m}{j}\mathcal{G}%
		_{n+j}^{(\alpha)}(x;\lambda )-
		\lambda (-1)^m \dfrac{n!m!}{(n+m+1)!}\mathcal{G}%
		_{0}^{(\alpha+1)}(x;\lambda ) 
		$$
	
\end{remark}

\section{ Limit case: $m$=$\protect\rho \,n,n\rightarrow
\infty.$}

In formulas (\ref{formule1teo2}, \ref{formule2teo2}, \ref{formule3teo2}) of Theorem \ref{teo2}, let us assume that $m=\rho\; n$ with $n$ going to infinity. After dividing (\ref{formule1teo2}) by $(n+m-p)!$, we get

$$\sum_{j=0}^{\max(m,n)} \left(\lambda \bin{n}{j}-(-1)^j\bin{m}{j}\right)\bin{n+m-j}{n+m-p}{\mathcal{B}%
	_{p-j}^{(\alpha+1)}(x;\lambda )}
=\sum_{j=0}^{m}  (-1)^j\bin{m}{j}\bin{n+m-j}{n+m-p}(p-j){\mathcal{B}%
	_{p-j-1}^{(\alpha)}(x;\lambda )},$$

 Making use of  $$\bin{n}{j}\underset{n\rightarrow \infty }{\sim }\frac{n^{j}}{j!},\;\;\;(j\geq 0)$$
and $$\bin{n+\rho \;n-j}{n+\rho\; n-p}\underset{n\rightarrow \infty }{\sim }\frac{(1+\rho
)^{p-j}}{(p-j)!}n^{p-j},\;\;\;( 0\leq j\leq p),$$
we obtain 

$$\sum_{j=0}^{\max(m,n)} \left(\lambda \dfrac{n^j}{j!}-(-\rho)^j\dfrac{n^j}{j!}\right)\frac{(1+\rho
	)^{p-j}}{(p-j)!}n^{p-j}{\mathcal{B}%
	_{p-j}^{(\alpha+1)}(x;\lambda )}
=\sum_{j=0}^{m}  (-\rho)^j\dfrac{n^j}{j!}\frac{(1+\rho
	)^{p-j}}{(p-j)!}n^{p-j}(p-j){\mathcal{B}%
	_{p-j-1}^{(\alpha)}(x;\lambda )},$$
after multiplying the two members by $\dfrac{p!}{n^p (1+\rho)^p}$
it arises

$$\sum_{j=0}^{p}\bin{p}{j}(\lambda -(-\rho )^{j})(1+\rho )^{-j}\mathcal{B}^{(\al+1)}%
_{-j}(x;\lambda ) =\!\!\sum_{j=0}^{p}p\bin{p-1}{j}(-\rho )^{j}(1+\rho )^{-j}\mathcal{B}^{(\al)}%
_{p-j-1}(x;\lambda ),
\label{RHOBERNOULLI}$$

Of course in the previous equation, $\rho$
is a rational number. But it is also valid if  $\rho$ is any complex number as proved in the following Theorem.

\begin{theorem}
\label{teo3}

For all complex number $\rho $ and for all integer $p$,  the following formulas hold:

\begin{eqnarray}
\sum_{j=0}^{p}\bin{p}{j}(\lambda -(-\rho )^{j})(1+\rho )^{-j}\mathcal{B}^{(\al+1)}%
_{p-j}(x;\lambda ) &=&\!\!\!\!\sum_{j=0}^{p}p\bin{p-1}{j}(-\rho )^{j}(1+\rho )^{-j}\mathcal{B}^{(\al)}%
_{p-j-1}(x;\lambda ) \notag\\
&=&p \;\mathcal{B}^{(\al)}%
_{p-1}\left(x-\dfrac{\rho}{1+\rho};\lambda \right),
\label{RHOBERNOULLI} \\
\sum_{j=0}^{p}\bin{p}{j}(\lambda +(-\rho )^{j})(1+\rho )^{-j}\mathcal{E}^{(\al+1)}%
_{p-j}(x;\lambda ) &=&2 \sum_{j=0}^{p}\bin{p}{j}(-\rho )^{j}(1+\rho )^{-j}\mathcal{E}^{(\al)}%
_{p-j}(x;\lambda ) \notag\\
&=&2 \;\mathcal{E}^{(\al)}%
_{p}\left(x-\dfrac{\rho}{1+\rho};\lambda \right),  \label{rhoeuler}\\
\sum_{j=0}^{p}\bin{p}{j}(\lambda +(-\rho )^{j}(1+\rho )^{-j}\mathcal{G}^{(\al+1)}%
_{p-j}(x;\lambda ) &=&\!\!\!\!2\;\sum_{j=0}^{p}p\bin{p-1}{j}(-\rho )^{j}(1+\rho )^{-j}\mathcal{G}^{(\al)}%
_{p-j-1}(x;\lambda )  \notag\\
&=&2 p \;\mathcal{G}^{(\al)}%
_{p-1}\left(x-\dfrac{\rho}{1+\rho};\lambda \right)\label{rhogenocchi}
\end{eqnarray}
\end{theorem}

\noindent  \textbf{Proof.}
It can be found that the generating functions of both sides of (\ref{RHOBERNOULLI}) are

 $\displaystyle \dfrac{(t(1+\rho))^{\al+1}}{(\lambda e^{t(1+\rho)}-1)^\al} e^{(x(1+\rho)-\rho)t}$,   of both sides of (\ref{rhoeuler}) are $\dfrac{2^{\al+1}}{\lambda e^{t(1+\rho)}+1)^\al} e^{(x(1+\rho)-\rho)t}$
 
  and 
of both sides of (\ref{rhogenocchi}) are
$\displaystyle \dfrac{(2t(1+\rho))^{\al+1}}{(\lambda e^{t(1+\rho)}+1)^\al} e^{(x(1+\rho)-\rho)t}$

\begin{remark}

 Set $X=x-\dfrac{\rho}{(1+\rho )},$\ the previous identities turn to 
\begin{equation*}
\lambda \mathcal{B}^{(\al+1)}_{p}(X+1;\lambda )-\mathcal{B}^{(\al+1)}_{p}(X;\lambda
)=p\mathcal{B}^{(\al)}_{p-1}(X;\lambda )(p\geq 0), \label{functionalrelationbernoulli}
\end{equation*}

\begin{equation*}
\lambda \mathcal{E}^{(\al+1)}_{p}(X+1;\lambda )+\mathcal{E}^{(\al+1)}_{p}(X;\lambda
)=2\mathcal{E}^{(\al)}_{p}(X;\lambda )(p\geq 0), \label{functionalrelationeuler}
\end{equation*}%
and
\begin{equation*}
\lambda \mathcal{G}^{(\al+1)}_{p}(X+1;\lambda )+\mathcal{G}^{(\al+1)}_{p}(X;\lambda
)=2p\mathcal{G}^{(\al)}_{p-1}(X;\lambda )(p\geq 0), \label{functionalrelationgenocchi}
\end{equation*}

which can be found in  \cite{apostol1951,luo2006,luosri2011} for the  first two relations.
\end{remark}

\section{Lacunary recurrence relation}

In this section we consider the particular case $\lambda=1$ and find lacunary relations for Bernoulli, Euler and Genocchi polynomials of higher order.

If $\rho=i \;(i^2=-1)$ then formulas (\ref{RHOBERNOULLI}), (\ref{rhoeuler}), (\ref{rhogenocchi})
reduce to the following  lacunary relations of length 4.

\begin{corollary}

\begin{eqnarray}
\sum_{k \equiv 2(4)}^p\bin{p}{k}2^{1-k/2}(-1)^{(k+2)/4}{B}_{p-k}^{(\al+1)}(x)&=&\Im \left(p\; B_{p-1}^{(\al)}\left(x-\dfrac{1+i}{2}\right)\right)\label{6.1}\\
\sum_{k \equiv 0(4)}^p\bin{p}{k}2^{1-k/2}(-1)^{k/4}{E}_{p-k}^{(\al+1)}(x)&=&\Re \left(2\; E_{p}^{(\al)}\left(x-\dfrac{1+i}{2}\right)\right)\\
\sum_{k \equiv 0(4)}^p\bin{p}{k}2^{1-k/2}(-1)^{k/4}{G}_{p-k}^{(\al+1)}(x)&=&\Re \left(2\;p\; G_{p-1}^{(\al)}\left(x-\dfrac{1+i}{2}\right)\right).
\end{eqnarray}

\end{corollary}

{\bf Proof
	}
Formula 	(\ref{RHOBERNOULLI}) with $\rho=i$ becomes

$$\sum_{j=0}^{p}\bin{p}{j}(1 -(-i )^{j})(1+i )^{-j}\mathcal{B}^{(\al+1)}%
_{p-j}(x ) =p \;\mathcal{B}^{(\al)}%
_{p-1}\left(x-\dfrac{i}{1+i} \right).$$

The coefficients $c_{j}:=(1 -(-i )^{j})(1+i )^{-j}=2^{-j/2}(e^{-i j \pi/4}-(-1)^j e^{i j \pi/4} ), $ satisfy:
\begin{eqnarray*}
	j &\equiv &0(4)\hspace{1cm}c_{j}=0, \\
	j &\equiv &1(4)\hspace{1cm}c_{j}=(-1)^{(j-1)/4}2^{(1-j)/2}, \\
	j &\equiv &2(4)\hspace{1cm}c_{j}=(-1)^{(j+2)/4}2^{(2-j)/2}\;i, \\
j &\equiv &3(4)\hspace{1cm}c_{j}=(-1)^{(j+1)/4}2^{(1-j)/2}.
\end{eqnarray*}

So, if we consider only the imaginary part of both sides of this equation, it provides (\ref{6.1}) which is a recurrence relation with a gap of length 4. 
The others formulas are proved    in the same manner.

\bibliographystyle{acm}

\bibliography{bibfile3}

\end{document}